# Identities for squared central binomial coefficients

**Khristo N. Boyadzhiev**
Department of Mathematics
Ohio Northern University, Ada, OH 45810
k-boyadzhiev@onu.edu

**Abstract** We prove four identities for the squared central binomial coefficients. The first three of them reflect certain transformation properties of the complete elliptic integrals of the first and the second kind, while the last one is based on properties of the Lagrange polynomials.

## 1. Introduction and main theorem

The central binomial coefficients are defined by

$$\binom{2n}{n} = \frac{(2n)!}{(n!)^2}, n = 0, 1, \ldots$$

There are various combinatorial results in the literature about these numbers, including interesting identities and congruences – see, for example, [2], [6, p.47], and [8]. They are usually obtained by combinatorial methods and special techniques, like in [6]. Riordan in his book [7, p. 130] proved the convolution type identity

$$\sum_{m=0}^{n} \binom{2m}{m} \frac{4^{n-m}}{1-2m} = \binom{2n}{n}.$$

It is interesting that similar identities exist for the squared central binomial coefficients. Here we prove several such identities, where (3) especially has almost the same structure as the one above. The proofs of the first three identities are based on certain transformation properties of the complete elliptic integrals of the first and second kind, $K(x)$ and $E(x)$ correspondingly. For instance, the first identity can be viewed as a numerical equivalent of Landen's transformation formula for $K(x)$. The proof of the third identity uses a transformation property for $E(x)$. The fourth identity is based on a property of the Legendre polynomials.

**Theorem**. For every nonnegative integer $n$ we have

$$\sum_{m=0}^{n} \binom{2m}{m}^2 \binom{n+m}{n-m} (-4)^{n-m} = \binom{2k}{k}^2 \text{ (for } n = 2k\text{)}, \tag{1}$$

and the sum is zero when $n$ is odd. Also,

$$\sum_{m=0}^{n} \binom{2m}{m}^2 \binom{n-1/2}{n-m} (-1)^m 16^{n-m} = \binom{2n}{n}^2. \tag{2}$$

$$\sum_{m=0}^{n} \binom{2m}{m}^2 \frac{16^{n-m}}{1-2m} = (2n+1)\binom{2n}{n}^2 \tag{3}$$

$$\sum_{m=0}^{2n} \binom{2m}{m}^2 \binom{2n+m}{2m} (-1)^m 4^{2n-m} = \binom{2n}{n}^2. \tag{4}$$

## 2. Proof of the theorem

*Proof.* For the proof of (1) and (2) we use the Taylor series expansion

$$K(x) = \frac{\pi}{2} \sum_{n=0}^{\infty} \binom{2n}{n}^2 \frac{x^{2n}}{4^{2n}}, \qquad (5)$$

of the complete elliptic integral of the first kind [1], [2, (3.127)], [4, (5.25.8)], [5].

$$K(x) = \int_0^{\frac{\pi}{2}} \frac{d\theta}{\sqrt{1 - x^2 \sin^2 \theta}} = \int_0^1 \frac{dt}{\sqrt{(1-t^2)(1-x^2 t^2)}}. \qquad (6)$$

When we expand $(1 - x^2 \sin^2 \theta)^{-1/2}$ in binomial series and integrate term by term, we obtain (5) by means of the well-known formula

$$\binom{-1/2}{n}(-1)^n = \binom{2n}{n}\frac{1}{4^n},$$

and the Wallis formula

$$\int_0^{\pi/2} \sin^{2n}\theta \, d\theta = \frac{\pi}{2^{2n+1}} \binom{2n}{n}.$$

The function $K(x)$ has the important transformation property (Landen's transformation [1], [5])

$$K(x) = \frac{1}{1+x} K\left(\frac{2\sqrt{x}}{1+x}\right)$$

which translates into the series identity

$$\sum_{n=0}^{\infty} \binom{2n}{n}^2 \frac{x^{2n}}{4^{2n}} = \frac{1}{1+x} \sum_{m=0}^{\infty} \binom{2m}{m}^2 \frac{x^m}{2^{2m}(1+x)^{2m}}.$$

With the substitution $x = 4t$ this becomes

$$\sum_{n=0}^{\infty} \binom{2n}{n}^2 t^{2n} = \sum_{m=0}^{\infty} \binom{2m}{m}^2 t^m (1+4t)^{-2m-1}. \qquad (7)$$

Expanding $(1+4t)^{-2m-1}$ we write

$$(1+4t)^{-2m-1} = \sum_{k=0}^{\infty} \binom{2m+k}{k}(-4)^k t^k$$

by using the property

$$\binom{-p}{k} = \binom{p+k-1}{k}(-1)^k.$$

Setting $n = m+k$ and exchanging the order of summation in the right hand side we obtain from (7)

$$\sum_{n=0}^{\infty} \binom{2n}{n}^2 t^{2n} = \sum_{n=0}^{\infty} t^n \left\{ \sum_{m=0}^{n} \binom{2m}{m}^2 \binom{n+m}{n-m} (-4)^{n-m} \right\}.$$

Identity (1) comes from here by comparing coefficients on both sides.

For the proof of (2) we use another transformation property of $K(x)$, namely,

$$K(ix) = \frac{1}{\sqrt{1+x^2}} K\left(\frac{x}{\sqrt{1+x^2}}\right) \tag{8}$$

This property is easy to verify directly in the integral (5). Applied to (4) it yields (with $x^2 = 16t$)

$$\sum_{n=0}^{\infty} \binom{2n}{n}^2 (-1)^n t^n = \sum_{m=0}^{\infty} \binom{2m}{m}^2 t^m (1+16t)^{-m-1/2}.$$

Again, expanding $(1+16t)^{-m-\frac{1}{2}}$ we write (with $m+k = n$)

$$\sum_{n=0}^{\infty} \binom{2n}{n}^2 (-1)^n t^n = \sum_{m=0}^{\infty} \binom{2m}{m}^2 t^m \left\{ \sum_{k=0}^{\infty} \binom{m+k-1/2}{k} (-1)^k 16^k t^k \right\}$$

$$= \sum_{n=0}^{\infty} t^n \left\{ \sum_{m=0}^{n} \binom{2m}{m}^2 \binom{n-1/2}{n-m} (-1)^{n-m} 16^{n-m} \right\}.$$

Obviously, (2) follows from here.

For the third identity (3) we use the complete elliptic integral of the second kind

$$E(x) = \int_{0}^{\frac{\pi}{2}} \sqrt{1-x^2 \sin^2 \theta}\; d\theta. \tag{9}$$

This function has the series expansion (see [1], [2, (3.128)])

$$E(x) = \frac{\pi}{2} \sum_{n=0}^{\infty} \binom{2n}{n}^2 \frac{x^{2n}}{4^{2n}(1-2n)}, \qquad (10)$$

and the property

$$\frac{\pi}{2} \sum_{n=0}^{\infty} \binom{2n}{n}^2 \frac{2n+1}{16^n} x^n = \frac{1}{1-x} E(\sqrt{x})$$

which is essentially entry (5.26.1) in Hansen's table [4]. From this

$$\sum_{n=0}^{\infty} \binom{2n}{n}^2 \frac{2n+1}{16^n} x^n = \frac{1}{1-x} \sum_{n=0}^{\infty} \binom{2n}{n}^2 \frac{x^n}{16^n(1-2n)}.$$

Now expanding $\frac{1}{1-x}$ in geometric power series and using Cauchy's rule for multiplication of two power series we obtain

$$\sum_{n=0}^{\infty} \binom{2n}{n}^2 \frac{2n+1}{16^n} x^n = \sum_{n=0}^{\infty} x^n \sum_{m=0}^{n} \binom{2m}{m}^2 \frac{1}{16^m(1-2m)}.$$

Identity (3) follows from here by comparing coefficients.

Next, let $P_n(x), n = 0, 1, \ldots$ be the Legendre polynomials [2], [3], [5], [8]. We have $P_n(-x) = (-1)^n P_n(x)$ and the representation [2, (3.135)]

$$P_n(x) = \sum_{k=0}^{n} \binom{n+k}{k}\binom{n}{k}\left(\frac{x-1}{2}\right)^k.$$

Using the simple identity

$$\binom{n+k}{k}\binom{n}{k} = \binom{n+k}{2k}\binom{2k}{k}$$

we can write

$$P_n(x) = \sum_{k=0}^{n} \binom{n+k}{2k}\binom{2k}{k}\left(\frac{x-1}{2}\right)^k.$$

Replacing here $x$ by $-\cos\theta$ yields

$$P_n(-\cos\theta) = \sum_{k=0}^{n} \binom{n+k}{2k}\binom{2k}{k}\left(\frac{1+\cos\theta}{2}\right)^k, \text{ or}$$

$$(-1)^n P_n(\cos\theta) = \sum_{k=0}^{n} \binom{n+k}{2k}\binom{2k}{k}(-1)^k \cos^{2k}\frac{\theta}{2}.$$

This equation we integrate from $0$ to $\pi$. Noticing that

$$\int_0^{\pi} \cos^{2k}\frac{\theta}{2}d\theta = 2\int_0^{\pi/2}\cos^{2k}\theta\, d\theta = \frac{\pi}{2^{2k}}\binom{2k}{k},$$

we arrive at the representation

$$(-1)^n \int_0^{\pi} P_n(\cos\theta)\, d\theta = \sum_{k=0}^{n}\binom{n+k}{2k}\binom{2k}{k}^2 \frac{(-1)^k}{4^k}.$$

The integral on the left hand side is nonzero only when $n$ is even and then

$$\int_0^{\pi} P_{2n}(\cos\theta)\, d\theta = \frac{\pi}{4^{2n}}\binom{2n}{n}^2$$

(see [3, 7.221(3) or [7.226(1)] ). Therefore,

$$\frac{\pi}{4^{2n}}\binom{2n}{n}^2 = \sum_{k=0}^{2n}\binom{2n+k}{2k}\binom{2k}{k}^2 \frac{(-1)^k}{4^k}$$

which is (4). □